# THE LARGEST COUNTABLE INDUCTIVE SET IS A MOUSE SET

MITCH RUDOMINER

ABSTRACT. Let $\kappa^{\mathbb{R}}$ be the least ordinal $\kappa$ such that $L_\kappa(\mathbb{R})$ is admissible. Let $A = \{x \in \mathbb{R} \mid (\exists \alpha < \kappa^{\mathbb{R}})$ s.t. $x$ is ordinal definable in $L_\alpha(\mathbb{R})\,\}$. It is well known that (assuming determinacy) $A$ is the largest countable inductive set of reals. Let $T$ be the theory: ZFC − Replacement + "There exists $\omega$ Woodin cardinals which are cofinal in the ordinals." $T$ has consistency strength weaker than that of the theory ZFC + "There exists $\omega$ Woodin cardinals", but stronger than that of the theory ZFC + "There exists $n$ Woodin Cardinals", for each $n \in \omega$. Let $\mathcal{M}$ be the canonical, minimal inner model for the theory $T$. In this paper we show that $A = \mathbb{R} \cap \mathcal{M}$. Since $\mathcal{M}$ is a *mouse*, we say that $A$ is a *mouse set*. As an application, we use our characterization of $A$ to give an inner-model-theoretic proof of Martin's theorem that $A$ is equal to the set of reals which are $\Sigma_n^*$ for some $n$.

## Contents

## 0. Introduction

This paper is best understood as being part of a certain programme in set theory—a programme whose goal is to explore the connections between descriptive set theory and inner model theory. We begin by describing some of the other contributions to this programme. The following several results are all of a common theme: From descriptive set theory one takes a countable, definable set of reals, $A$. It is then shown that $A = \mathbb{R} \cap \mathcal{M}$, where $\mathcal{M}$ is some canonical model from inner model theory. In technical terms, $\mathcal{M}$ is a mouse. Consequently we say that $A$ is a *mouse set*.

We begin with some results due to D. A. Martin, J. R. Steel, and H. Woodin. First a definition. If $\xi$ is a countable ordinal and $x \in \mathbb{R}$ and $n \geq 2$, then we say that $x \in \Delta^1_n(\xi)$ if for every real $w \in \text{WO}$ such that $|w| = \xi$, $x \in \Delta^1_n(w)$. Set

$$A_n = \left\{ x \in \mathbb{R} \,\middle|\, (\exists \xi < \omega_1) \, x \in \Delta^1_n(\xi) \right\}.$$

In short, $A_n$ is the set of reals which are $\Delta^1_n$ in a countable ordinal. Under the hypothesis of Projective Determinacy (PD), the sets $A_n$ are of much interest to descriptive set theorists. If $n \geq 2$ is even then $A_n = C_n$, the largest countable $\Sigma^1_n$ set. If $n \geq 3$ is odd then $A_n = Q_n$, the largest countable $\Pi^1_n$ set which is closed downward under $\Delta^1_n$ degrees. The sets $A_n$ have been studied extensively by descriptive set theorists. See for example [Kech] and [KeMaSo].

The work of Martin, Steel, and Woodin shows that there is a connection between the sets $A_n$, and certain inner models of the form $L[\vec{E}]$, where $\vec{E}$ is a sequence of extenders. In the paper "Projectively Wellordered Inner Models" [St3], Martin, Steel, and Woodin prove the following theorem:

**Theorem 0.1** (Martin, Steel, Woodin). *Let $n \geq 1$ and suppose that there are $n$ Woodin cardinals with a measurable cardinal above them. Let $\mathcal{M}_n$ be the canonical $L[\vec{E}]$ model with $n$ Woodin cardinals. Then*

$$\mathbb{R} \cap \mathcal{M}_n = A_{n+2}.$$

The above theorem is also true with $n = 0$. The "canonical $L[\vec{E}]$ model with 0 Woodin cardinals" is just the model $L$. If there is a measurable cardinal, then $A_2 = C_2 = \mathbb{R} \cap L$. Actually, this result is true under the weaker hypothesis that $\mathbb{R} \cap L$ is countable.

Now let $A^*$ be the collection of reals which are ordinal definable in $L(\mathbb{R})$. To emphasize the parallel with the sets $A_n$ above, note that assuming that every game in $L(\mathbb{R})$ is determined ($\text{AD}^{L(\mathbb{R})}$), we have:

$$A^* = \left\{\, x \in \mathbb{R} \;\middle|\; (\exists \xi < \omega_1)\, x \in (\Delta^2_1)^{L(\mathbb{R})}(\xi) \,\right\},$$

and $A^*$ is the largest countable $(\Sigma^2_1)^{L(\mathbb{R})}$ set of reals.

In the paper "Inner Models with Many Woodin Cardinals" [St2], Steel and Woodin prove the following theorem:

**Theorem 0.2** (Steel, Woodin). *Suppose that there are $\omega$ Woodin cardinals with a measurable cardinal above them. Let $\mathcal{M}_\omega$ be the canonical $L[\vec{E}]$ model with $\omega$ Woodin cardinals. Then*

$$\mathbb{R} \cap \mathcal{M}_\omega = A^*.$$

Theorems 0.1 and 0.2 are obviously similar to each other. Let us establish some terminology which will allow us to describe this similarity. Given a pointclass $\Delta$ let us set:

$$A_\Delta = \left\{\, x \in \mathbb{R} \;\middle|\; (\exists \xi < \omega_1)\, x \in \Delta(\xi) \,\right\}.$$

Theorems 0.1 and 0.2 both establish results of the form: $A_\Delta$ is a mouse set. Theorem 0.1 does this for $\Delta = \Delta^1_n$. Theorem 0.2 does this for $\Delta = (\Delta^2_1)^{L(\mathbb{R})}$. Furthermore, the two theorems do not merely show that there *exists some* mouse $\mathcal{M}$ such that $A_\Delta = \mathbb{R} \cap \mathcal{M}$. The theorems actually describe $\mathcal{M}$ in terms of the large cardinal axioms that it satisfies.

There is quite a bit of room between theorems 0.1 and 0.2. One way to describe this room is in terms of the levels of the model $L(\mathbb{R})$. For convenience, we will use the Jensen $J$-hierarchy for $L(\mathbb{R})$. We define $J_1(\mathbb{R})$ to be $\mathbb{R} \cup \text{HF}$. For $\lambda$ a limit ordinal we set $J_\lambda(\mathbb{R}) = \bigcup_{\alpha < \lambda} J_\alpha(\mathbb{R})$. For $\alpha \geq 1$ we set $J_{\alpha+1}(\mathbb{R})$ to be the rudimentary closure of $J_\alpha(\mathbb{R}) \cup \{J_\alpha(\mathbb{R})\}$. Then $L(\mathbb{R}) = \bigcup_\alpha J_\alpha(\mathbb{R})$. See [St1] or [Ru] for a more detailed discussion. Now let $\delta = (\delta^2_1)^{L(\mathbb{R})}$ be the least ordinal such that $J_\delta(\mathbb{R})$ is a $\Sigma_1$ elementary submodel of $L(\mathbb{R})$. Then we can characterize the room between theorems 0.1 and 0.2 by pointing out that the first theorem is concerned with the pointclasses $\Delta_n(J_1(\mathbb{R}))$, and the second theorem is concerned with the pointclass $\Delta_1(J_\delta(\mathbb{R}))$. Assuming $\text{AD}^{L(\mathbb{R})}$, $\delta$ is a large cardinal in $L(\mathbb{R})$ (for example it is an inaccessible limit of inaccessibles.) From this point of view there is quite a bit of room between theorems 0.1 and 0.2.

Given the above discussion, it is natural to conjecture that $A_\Delta$ is a mouse set for $\Delta = \Delta_n(J_\alpha(\mathbb{R}))$, for *all* $\alpha$ and *all* $n$. In [Ru] we do make such a conjecture (in a slightly more precise formulation.) In that paper we are not able to fully prove the conjecture, but we are able to prove the conjecture for *some* $\alpha$ and *some* $n$. Our main theorem in [Ru] is similar to theorems 0.1 and 0.2 above in that we not only show that there *exists some* mouse $\mathcal{M}$ such that $A_\Delta = \mathbb{R} \cap \mathcal{M}$, we actually describe $\mathcal{M}$ in terms of the large cardinal axioms which it satisfies. Independently, and using entirely different techniques, Woodin has shown that $A_\Delta$ is a mouse set for $\Delta = \Delta_1(J_\lambda(\mathbb{R}))$, for all limit ordinals $\lambda$. See [St4]. Woodin's proof shows only that there *exists some* mouse $\mathcal{M}$ such that $A_\Delta = \mathbb{R} \cap \mathcal{M}$. The mouse $\mathcal{M}$ is not described in terms of the large cardinal axioms which it satisfies.

Now we turn to the topic of the current paper. In this paper we will show that $A_\Delta$ is a mouse set, for $\Delta = $ the pointclass of *hyperprojective* sets. (See below for the definition of hyperprojective.) Assuming determinacy, $A_\Delta$ is the largest countable inductive set. Thus there are obvious parallels between our result, and the results of theorems 0.1 and 0.2 above. We may also compare our result with those of theorems 0.1 and 0.2 by looking at levels of $L(\mathbb{R})$. Another way of describing the main result of this paper is to say that we will show that $A_\Delta$ is a mouse set, for $\Delta = \Delta_1(J_\kappa(\mathbb{R}))$, where $\kappa = \kappa^\mathbb{R}$ is the least ordinal such that $J_\kappa(\mathbb{R})$ is admissible. Since $1 < \kappa^\mathbb{R} < (\delta_1^2)^{L(\mathbb{R})}$, the result of this paper is strictly *between* those of theorem 0.1 and and theorem 0.2 above.

We will feel free to use terms and concepts from inner model theory. In particular, we will expect the reader to have some familiarity with the papers [MiSt] and [St2].

Some of the research for this paper was done while I was a graduate student at UCLA, although the results do not appear in my PhD thesis. I would like to thank my thesis advisor, Professor John Steel.

### The Pointclass of Inductive Sets

The pointclass of inductive sets has several different equivalent descriptions. (Throughout this paper we shall only be using the term "inductive" in the *lightface* sense. This is what is called "absolutely inductive" in [Mo].) In section 7C of [Mo] the inductive sets are defined as the *positive analytical inductive sets on* $\mathbb{R}$. It is shown there that the class of all inductive sets is the smallest Spector pointclass which is closed under both $\forall^\mathbb{R}$ and $\exists^\mathbb{R}$. The compliment of an inductive set is called *co-inductive.* If $X$ is both inductive and co-inductive then we say that $X$ is *hyperprojective.* As we shall not have any need for the concept of "positive analytical inductive," we shall not go into any more detail about it here. Instead we give two other equivalent characterizations of the inductive sets.

**Definition 0.3.** $\kappa^\mathbb{R}$ is the least ordinal $\kappa$ so that $J_\kappa(\mathbb{R})$ is admissible.

For convenience we make the following convention.

**Notation 0.4.** For the rest of the paper $\kappa$ will refer exclusively to $\kappa^\mathbb{R}$.

A set $X \subset \mathbb{R}$ is inductive iff $X$ is $\Sigma_1$ definable over the model $J_\kappa(\mathbb{R})$, with the point $\mathbb{R}$ as a parameter. As we shall always allow the point $\mathbb{R}$ as a parameter in definitions over the models $J_\alpha(\mathbb{R})$, we suppress mention of it in our notation. Thus we write that $X$ is inductive iff $X$ is $\Sigma_1(J_\kappa(\mathbb{R}))$. (See [St1] or [Ru] for a more detailed discussion of definability over the models $J_\alpha(\mathbb{R})$.) Similarly, $X$ is co-inductive iff $X$ is $\Pi_1(J_\kappa(\mathbb{R}))$ and $x$ is hyperprojective iff $X$ is $\Delta_1(J_\kappa(\mathbb{R}))$.

**Definition 0.5.** $A_{\mathrm{HYP}}$ = the set of all $x \in \mathbb{R}$ such that $(\exists \alpha < \kappa^\mathbb{R})$ so that $x$ is definable over $J_\alpha(\mathbb{R})$ from ordinal parameters.

It is not difficult to see that $x \in A_{\mathrm{HYP}}$ iff $\{x\}$ is $\Sigma_1(J_\kappa(\mathbb{R}), \{\xi\})$, for some ordinal $\xi$, iff $\{x\}$ is $\Delta_1(J_\kappa(\mathbb{R}), \{\xi\})$. Assuming determinacy $A_{\mathrm{HYP}}$ is countable, and so we may take $\xi$ to be a countable ordinal. Thus $A_{\mathrm{HYP}}$ is what we were above calling $A_\Delta$, for $\Delta$ = the pointclass of hyperprojective sets. Also assuming determinacy, it is not difficult to see that $A_{\mathrm{HYP}}$ is the largest countable inductive set.

In this paper we will show that $A_{\mathrm{HYP}}$ is a mouse set, that is that $A_{\mathrm{HYP}} = \mathbb{R} \cap \mathcal{M}$ for some mouse $\mathcal{M}$. There are two halves to this proof. In section 1 we use a comparison argument to show that $\mathbb{R} \cap \mathcal{M} \subseteq A_{\mathrm{HYP}}$, and in section 2 we prove a correctness theorem which implies that $A_{\mathrm{HYP}} \subseteq \mathbb{R} \cap \mathcal{M}$.

Let $\Sigma_0^* =$ the class of unions of inductive and coinductive sets. Let $\Pi_n^*$ be the class of complements of $\Sigma_n^*$ sets and let $\Sigma_{n+1}^*$ be the class of projections of $\Pi_n^*$ sets. Equivalently, $\Sigma_n^* = \Sigma_{n+1}(J_\kappa(\mathbb{R}))$, for $n \geq 1$. (To see this use the fact that $J_\kappa(\mathbb{R})$ projects to $\mathbb{R}$. See [St1].) In [Ma2], Martin shows that $A_{\text{HYP}}$ is equal to the set of reals which are $\Sigma_n^*$ for some $n$. In section 3 we reprove one direction of Martin's theorem using purely inner-model-theoretic techniques. Using our characterization of $A_{\text{HYP}}$ as $\mathbb{R} \cap \mathcal{M}$ for some premouse $\mathcal{M}$, we show that every $\Sigma_n^*$ real $x$ is in $A_{\text{HYP}}$ by showing that $x \in \mathcal{M}$.

We shall need one other characterization of the inductive sets: A set $A \subseteq \mathbb{R}$ is inductive iff $A$ is $\supset^\mathbb{R}$-open. That is iff there is a function $x \mapsto G_x$ such that for all reals $x$, $x \in A \leftrightarrow$ player $I$ wins $G_x$, where $G_x$ is an open game on $\mathbb{R}$ which is continuously associated to $x$. More precisely, $A$ is inductive iff there is some arithmetic set $P \subseteq \omega^{<\omega} \times (\omega^{<\omega})^{<\omega}$ such that $x \in A$ iff

$$(\exists y_1 \in \mathbb{R})(\forall y_2 \in \mathbb{R})(\exists y_3 \in \mathbb{R}) \cdots (\exists n \in \omega) \; P(x \restriction n, y_1 \restriction n, y_2 \restriction n, \ldots, y_n \restriction n).$$

Similarly, $A$ is co-inductive iff $A$ is $\supset^\mathbb{R}$-closed.

**Fairly Small Premice**

In this paper we will show that $A_{\text{HYP}}$ is a mouse set. We shall not only show that there exists some mouse $\mathcal{M}$ such that $A_{\text{HYP}} = \mathbb{R} \cap \mathcal{M}$, we shall describe $\mathcal{M}$ in terms of the large cardinal axioms which it satisfies. By examining theorems 0.1 and 0.2 above, we can see that $\mathcal{M}$ must satisfy a large cardinal axiom which is stronger than ZFC + "there exists $n$ Woodin cardinals", for each $n$, but weaker than ZFC + "there exists $\omega$ Woodin cardinals". To obtain the large cardinal axiom which $\mathcal{M}$ will satisfy, we start with the hypothesis ZFC + "there exists $\omega$ Woodin cardinals", and we remove the assumption that there are any ordinals above the supremum of the $\omega$ Woodin cardinals.

The mice used in section 4 of [St2] are called $\omega$-small. The mice used in [St3] are called $n$-small. We need a notion which fits between these two notions. Without any thought as to its use beyond this paper, we will adopt the term *fairly small*.

**Definition 0.6.** Let $\mathcal{M}$ be a premouse. Then we will say that $\mathcal{M}$ is *fairly big* iff there is an initial segment $\mathcal{N} \trianglelefteq \mathcal{M}$ such that there is an increasing $\omega$-sequence of ordinals $\delta_1 < \delta_2 < \delta_3, \ldots$ such that each $\delta_i$ is a Woodin cardinal of $\mathcal{N}$. If $\mathcal{M}$ is not fairly big then we will say the $\mathcal{M}$ is *fairly small*.

## 1. Comparison Lemma for $\partial^{\mathbb{R}}$-Closed Iterable Mice

In this section we will give half of the proof that $A_{\text{HYP}}$ is a mouse set. We will show that if $x \in \mathbb{R} \cap \mathcal{M}$, where $\mathcal{M}$ is a fairly small, iterable premouse, then $x \in A_{\text{HYP}}$. This is true because, essentially, $\mathcal{N} \in A_{\text{HYP}}$, where $\mathcal{N} \trianglelefteq \mathcal{M}$ is the $\trianglelefteq$-least initial segment of $\mathcal{M}$ containing $x$. The reason that $\mathcal{N}$ is simply definable is that $\mathcal{N}$ is the unique premouse of its ordinal height which is sound, projects to $\omega$, and is iterable. This would be enough to see that (the real coding) $\mathcal{N}$ was in $A_{\text{HYP}}$, if we could see that the property of being iterable was sufficiently simply definable. As in section 4 of [St2] and in [St3], we will define a weakening of the notion of full iterability which is simply definable, and yet is strong enough to compare fairly small premice. Our weakened notion of iterability is best described using the notion of an *iteration game*.

In section 1 of [St2], Steel defines $W\mathcal{G}_n(\mathcal{M}, \theta)$, the $n$-maximal, weak iteration game on $\mathcal{M}$, of length $\theta$. For our purposes here, we will need a slight modification of the game $W\mathcal{G}_0(\mathcal{M}, \omega)$. The set $\{\mathcal{M} \mid II \text{ wins } W\mathcal{G}_0(\mathcal{M}, \omega)\}$ is $\partial^{\mathbb{R}}$-$\Pi_1^1$. We need a notion of iterability which is $\partial^{\mathbb{R}}$-closed. Below we define an iteration game called the *closed iteration game on* $\mathcal{M}$.

Let $\mathcal{M}$ be a countable, premouse. The *closed iteration game on* $\mathcal{M}$ is identical to $W\mathcal{G}_0(\mathcal{M}, \omega)$, except that if neither player loses during the first $\omega$ rounds of play, then we automatically declare that player $II$ has won. That is, we do not demand that the direct limit model be defined and wellfounded in order for player $II$ to win. Thus the game is essentially a game on $\mathbb{R}$ with closed payoff for player $II$. Here are the rules of the closed iteration game in more detail: The game is played in $\omega$ rounds. Before beginning round $n < \omega$ we have a premouse $\mathcal{M}_n$, and an integer $k_n$ such that $\mathcal{M}_n$ is $k_n$-sound. We get started by setting $\mathcal{M}_0 = \mathcal{M}$ and $k_0 = 0$. Round $n$ is played as follows. Player $I$ begins by playing a countable, putative $k_n$-maximal iteration tree $\mathcal{T}$ on $\mathcal{M}_n$. Player $II$ can then either accept $\mathcal{T}$ or play a maximal well-founded branch $b$ of $\mathcal{T}$, with the proviso that he cannot accept $\mathcal{T}$ if it has a last, ill-founded model. If it is impossible for $II$ to play then he loses at round $n$. Suppose that $II$ does not lose at round $n$. If $II$ acceptss $\mathcal{T}$, then we let $\mathcal{M}_{n+1}$ be the last model of $\mathcal{T}$. If $II$ picks a maximal wellfounded branch $b$ of $\mathcal{T}$, then we set $\mathcal{M}_{n+1} = \mathcal{M}_b^{\mathcal{T}}$. In either case we let $k_{n+1}$ be the degree of $\mathcal{M}_{n+1}$ in $\mathcal{T}$. If $II$ does not lose at any round $n < \omega$ then $II$ wins.

**Definition 1.1.** Let $\mathcal{M}$ be a countable premouse. Then $\mathcal{M}$ is $\partial^{\mathbb{R}}$-*closed iterable* iff player $II$ wins the closed iteration game on $\mathcal{M}$.

**Remark 1.2.** The set of (reals coding) premice $\mathcal{M}$ such that player $II$ wins the closed iteration game on $\mathcal{M}$ is a a $\partial^{\mathbb{R}}$-closed set. In other words the set is coinductive. In other words the set is $\Pi_1(J_\kappa(\mathbb{R}))$.

Let $\mathcal{M}$ and $\mathcal{N}$ be countable premice. We say that $\mathcal{M}$ and $\mathcal{N}$ *can be compared* iff there is a countable iteration tree $\mathcal{T}$ on $\mathcal{M}$ of length $\theta + 1$ and a countable iteration tree $\mathcal{U}$ on $\mathcal{N}$ of length $\mu + 1$ such that either

(1) $\mathcal{M}_\theta^{\mathcal{T}}$ is an initial segment of $\mathcal{M}_\mu^{\mathcal{U}}$, and $D^{\mathcal{T}} \cap [0, \theta]_{\mathcal{T}} = \emptyset$, or

(2) $\mathcal{N}_\mu^{\mathcal{U}}$ is an initial segment of $\mathcal{N}_\theta^{\mathcal{T}}$, and $D^{\mathcal{U}} \cap [0, \mu]_{\mathcal{U}} = \emptyset$.

The main result in this section says that if $\mathcal{M}$ is fully iterable and fairly small, and $\mathcal{N}$ is $\partial^{\mathbb{R}}$-closed iterable, then $\mathcal{M}$ and $\mathcal{N}$ can be compared. To see that our comparison process terminates, we will need a certain amount of *generic absoluteness*. We will get the generic absoluteness we need by assuming that there exists infinitely many Woodin cardinals in $V$, and then using the machinery of Woodin's *stationary tower forcing*. See Chapter 9 of [Ma1] for a thorough treatment of this machinery. Below we simply quote one fact from the theory of stationary tower forcing.

**Proposition 1.3** (Woodin). *Let $\langle \delta_n \mid n \in \omega \rangle$ be a strictly increasing sequence of Woodin cardinals. Let $\lambda = \sup_n \delta_n$. Let $\gamma > \lambda$ be any ordinal. Then there is a generic extension* *of the universe, $V[G]$, such that in $V[G]$ there is an elementary embedding $j : V \to$ Ult such that*

(i) *$\gamma$ is in the wellfounded part of* Ult,

(ii) *$\mathrm{crit}(j) = \omega_1^V$, and $j(\omega_1^V) = \lambda$, and*

(iii) *letting $\mathbb{R}^*$ be the reals of* Ult, *we have that $\mathbb{R}^*$ is also the set of reals in a symmetric collapse of $V$ up to $\lambda$.*

See section 9.5 of [Ma1] for a proof of this proposition. Now we turn to the main result in this section.

**Remark 1.4.** Below we will use the notion of *realizability*. A mouse $\mathcal{M}$ is *realizable* if, roughly speaking, $\mathcal{M}$ can be embedded into a mouse $\bar{\mathcal{M}}$, where $\bar{\mathcal{M}}$ has the property that the extenders on the $\bar{\mathcal{M}}$ sequence all come from full extenders in $V$. See [St2] for the precise definition of realizable. In [St2] it is shown that if $\mathcal{M}$ is realizable, then $\mathcal{M}$ is "sufficiently iterable." In particular we will use the fact that if $\mathcal{M}$ is realizable then $\mathcal{M}$ is $\partial^{\mathbb{R}}$-closed iterable.

**Lemma 1.5** (Comparison Lemma). *Assume that there exists $\omega$ Woodin cardinals. Let $\mathcal{M}$ be countable, fairly small, realizable premouse. Let $\mathcal{N}$ be a countable, $\partial^{\mathbb{R}}$-closed iterable premouse. Then $\mathcal{M}$ and $\mathcal{N}$ can be compared.*

*Proof.* The proof of Theorem 1.10 from [St2] shows this. We will give a sketch of the details. In the language of that proof, when comparing $\mathcal{M}$ with $\mathcal{N}$, when more than one cofinal realizable branch appears on the $\mathcal{M}$-side, or more than one cofinal $\partial^{\mathbb{R}}$-closed iterable branch appears on the $\mathcal{N}$-side, we dovetail in a new comparison which guarantees that a certain ordinal $\delta$ is Woodin in the common "lined up part" of the family of models we are comparing. Thus every time we dovetail in a new comparison, we get another Woodin cardinal. Suppose we had to dovetail in a new comparison infinitely many times. Since $\mathcal{M}$ is realizable, on the $\mathcal{M}$ side we have a wellfounded direct limit model $\mathcal{M}'$. Because we had to dovetail in a new comparison infinitely many times, in $\mathcal{M}'$ there are infinitely many ordinals $\delta$ which are Woodin. Thus $\mathcal{M}'$ is not fairly small. But this contradicts our assumption that $\mathcal{M}$ is fairly small. Thus we did not have to dovetail in a new comparison infinitely many times.

In the proof of Theorem 1.10 from [St2], we would then be in Case 2 of the proof of the Claim. That is, we still need to see that our comparison does not last for $\omega_1$ stages. Suppose that our comparison does last for $\omega_1$ stages. Let $\mathcal{T}$ on $\mathcal{M}$ and $\mathcal{U}$ on $\mathcal{N}$ be the trees of length $\omega_1$ which are generated by the comparison. We can derive the usual contradiction if we can show that there is a cofinal branch $b$ of $\mathcal{T}$ and a cofinal branch $c$ of $\mathcal{U}$. Since $\mathcal{M}$ is realizable, the proof of Theorem 1.10 shows that there is a cofinal branch $b$ of $\mathcal{T}$. To see that there is a cofinal branch $c$ of $\mathcal{U}$, we will need to use our hypothesis that there exists $\omega$ Woodin cardinals in $V$. Let $\lambda$ be the supremum of the $\omega$ Woodin cardinals.

Let $\mathbb{R}'$ be the set of reals in some symmetric collapse of $V$ up to $\lambda$. Then there is some ordinal $\alpha$ such that $J_\alpha(\mathbb{R}')$ is admissible. Let $\alpha_0$ be the least such ordinal. Since the collapse forcing is homogeneous, $\alpha_0$ does not depend on our particular choice of $\mathbb{R}'$.

Let $\gamma > \alpha_0$ be any ordinal. Recall that we have set $\kappa = \kappa^{\mathbb{R}}$. Let $j : V \to \text{Ult}$ be as in Proposition 1.3 above. Let $\mathbb{R}^*$ be the reals of Ult. By Proposition 1.3, $\mathbb{R}^*$ is also the set of reals in a symmetric collapse of $V$ up to $\lambda$. Thus $\alpha_0$ is the least $\alpha$ such that $J_\alpha(\mathbb{R}^*)$ is admissible. Now $\alpha_0$ is in the wellfounded part of Ult. So $\text{Ult} \models$ "$\alpha_0$ is the least $\alpha$ such that $J_\alpha(\mathbb{R})$ is admissible." Thus $j(\kappa) = \alpha_0$.

Now $j(\mathcal{U})$ is an iteration tree on $\mathcal{N}$ which properly extends $\mathcal{U}$. By case hypothesis, there is some ordinal $\beta < \omega_1$ such that for all $\eta$ with $\beta < \eta < \omega_1^V$, $\mathcal{U} \upharpoonright \eta$ has a unique, cofinal, $\partial^{\mathbb{R}}$-closed iterable branch. (Otherwise we would have had to dovetail in a new comparison $\omega_1$ times.) Since $\partial^{\mathbb{R}}$-closed $= \Pi_1(J_\kappa(\mathbb{R}))$ and $j : J_\kappa(\mathbb{R}) \to J_{\alpha_0}(\mathbb{R}^*)$ is fully elementary, we have that for all $\eta$ with $\beta < \eta < \omega_1^{J_{\alpha_0}(\mathbb{R}^*)}$, $j(\mathcal{U}) \upharpoonright \eta$ has a unique, cofinal branch that is $\partial^{\mathbb{R}}$-closed iterable in $J_{\alpha_0}(\mathbb{R}^*)$. Since $j(\mathcal{U}) \upharpoonright \omega_1^V = \mathcal{U}$, there is a cofinal branch $c$ of $\mathcal{U}$ which is definable in $L(\mathbb{R}^*)$ from $\mathcal{U}$ and $\alpha_0$. As the collapse forcing is homogeneous, $c \in V$. $\square$

**Corollary 1.6.** *Assume that there exists $\omega$ Woodin cardinals. Let $\mathcal{M}$ be a countable, realizable, fairly small premouse. Then every real in $\mathcal{M}$ is ordinal definable over $J_\beta(\mathbb{R})$, for some $\beta < \kappa^{\mathbb{R}}$.*

*Proof.* Fix a real $x$ in $\mathcal{M}$. Let $\gamma$ be the rank of $x$ in the order of construction of $\mathcal{M}$. We will show that $x$ is definable from $\gamma$. Notice that $x$ is the unique real $x'$ such that:

there exists a countable, fairly small, $\supseteq^{\mathbb{R}}$-closed iterable premouse $\mathcal{N}$ such that

$x' \in \mathcal{N}$ and $x'$ is the $\gamma$th real in the order of construction of $\mathcal{N}$.

(proof: $\mathcal{M}$ witnesses that the statement above is true of the real $x$. Suppose that $\mathcal{N}$ witnessed that the statement above was true of some other real $x'$. By our Comparison Lemma, $\mathcal{M}$ and $\mathcal{N}$ can be compared. This implies that $x' = x$.) Since $\Pi_1(J_\kappa(\mathbb{R}))$ is closed under real quantification, we have that $\{x\}$ is $\Pi_1(J_\kappa(\mathbb{R}), \gamma)$. Fix a $\Sigma_1$ formula $\varphi$ such that $x$ is the unique real $x'$ so that $J_\kappa(\mathbb{R}) \models \neg\varphi[x', \gamma, \mathbb{R}]$. Let $f : \mathbb{R} \to \kappa^{\mathbb{R}}$ be defined by: $f(x) = 0$, and for $y \neq x$, $f(y) =$ the least $\alpha < \kappa^{\mathbb{R}}$ such that $J_\alpha(\mathbb{R}) \models \varphi[y, \gamma, \mathbb{R}]$. Notice that $f$ is $\Sigma_1(J_\kappa(\mathbb{R}), \{\gamma, x\})$. Thus the range of $f$ is not cofinal in $\kappa^{\mathbb{R}}$. Let $\beta = \sup(\text{ran}(f))$. Then $x$ is the unique real $x'$ such that $J_\beta(\mathbb{R}) \models \neg\varphi[x', \gamma, \mathbb{R}]$. Thus $x$ is ordinal definable over $J_\beta(\mathbb{R})$. □

## 2. INDUCTIVE CORRECTNESS

In this section we will complete the proof that $A_{\text{HYP}}$ is a mouse set. We will show that if $\mathcal{M}$ is iterable and fairly big, then $A_{\text{HYP}} \subseteq \mathbb{R} \cap \mathcal{M}$. We will need to use Woodin's stationary tower forcing. See chapter 9 of [Ma1] for a thorough treatment of stationary tower forcing.

**Definition 2.1.** Let $\delta$ be a limit ordinal. Then $\mathbb{Q}_\delta$ is the "$\mathbb{Q}$-version" of the stationary tower forcing up to $\delta$. More specifically, $x \in \mathbb{Q}_\delta$ iff $x \in V_\delta$ and $x$ is a stationary set and every element of $x$ is countable. The ordering on $\mathbb{Q}_\delta$ is described in chapter 9 of [Ma1].

We will need to use the following technical fact about the $\mathbb{Q}_\delta$.

**Lemma 2.2** (Woodin). *Let $\delta$ be a Woodin cardinal, and let $\rho > \delta$ be a limit ordinal. There is a condition $\hat{S} \in \mathbb{Q}_\rho$ such that for all $T \in \mathbb{Q}_\delta$, $\hat{S}$ is compatible with $T$ in $\mathbb{Q}_\rho$. Furthermore, if $G$ is $V$-generic over $\mathbb{Q}_\rho$ with $\hat{S} \in G$, then $G \cap \mathbb{Q}_\delta$ is $V$-generic over $\mathbb{Q}_\delta$.*

*Proof.* See chapter 9 of [Ma1]. $\square$

Using the above-mentioned fact, we will now prove another technical lemma about the $\mathbb{Q}_\delta$. Suppose that $\delta_1 < \delta_2$ are Woodin cardinals. The next lemma says, roughly, that a $\mathbb{Q}_{\delta_1}$-generic object can be extended to a $\mathbb{Q}_{\delta_2}$-generic object in such a way as to absorb some other given generic object.

**Lemma 2.3.** *Let $\mathfrak{M}$ be a transitive inner model of some sufficiently large fragment of ZFC. Let $\delta_1 < \delta_2$ be countable ordinals such that in $\mathfrak{M}$, $\delta_1$ and $\delta_2$ are both Woodin cardinals. Let $\mathbb{Q}_1 = \mathbb{Q}_{\delta_1}^{\mathfrak{M}}$ and $\mathbb{Q}_2 = \mathbb{Q}_{\delta_2}^{\mathfrak{M}}$. Let $\mathbb{P} \in \mathfrak{M}$ be any partial order such that $\mathrm{card}(\mathbb{P})^{\mathfrak{M}} < \delta_2$. Let $G_1$ and $H$ be sets such that $G_1$ is $\mathfrak{M}$-generic over $\mathbb{Q}_1$ and $H$ is $\mathfrak{M}[G_1]$-generic over $\mathbb{P}$. Then there exists a $G_2$ such that*

- *$G_2$ is $\mathfrak{M}$-generic over $\mathbb{Q}_2$.*
- *$G_2 \cap \mathbb{Q}_1 = G_1$.*
- *$H \in \mathfrak{M}[G_2]$.*

*Proof.* By Lemma 2.2 above, there is a condition $\hat{S} \in \mathbb{Q}_2$ such that for all $T \in \mathbb{Q}_1$, $\hat{S}$ is compatible with $T$ in $\mathbb{Q}_2$. Furthermore, if $G$ is $\mathfrak{M}$-generic over $\mathbb{Q}_2$ with $\hat{S} \in G$, then $G \cap \mathbb{Q}_1$ is $\mathfrak{M}$-generic over $\mathbb{Q}_1$. Fix such a condition $\hat{S}$. Let $\mathbb{Q}_2'$ be the collection of elements $T \in \mathbb{Q}_2$ such that $T$ is compatible with $\hat{S}$. So we have that $\mathbb{Q}_1 \subset \mathbb{Q}_2' \subset \mathbb{Q}_2$. Below we will use the notion of a *complete* embedding from one partial order into another. (See Definition 7.1 on page 218 of [Kunen].)

**Claim 1.** $\mathbb{Q}_1$ *is a complete suborder of* $\mathbb{Q}_2'$.

*Proof of Claim 1.* Let $G$ be $\mathfrak{M}$-generic over $\mathbb{Q}_2'$. Then $G$ is $\mathfrak{M}$-generic over $\mathbb{Q}_2$, and $\hat{S} \in G$. So $G \cap \mathbb{Q}_1$ is $\mathfrak{M}$-generic over $\mathbb{Q}_1$. In summary we have that whenever $G$ is generic over $\mathbb{Q}_2'$, $G \cap \mathbb{Q}_1$ is generic over $\mathbb{Q}_1$. It is an easy exercise to see that this condition is equivalent to the fact that $\mathbb{Q}_1$ is a complete suborder of $\mathbb{Q}_2'$. $\diamondsuit$

Now fix $G_1$ and $H$ so that $G_1$ is $M$-generic over $\mathbb{Q}_1$ and $H$ is $M[G_1]$-generic over $\mathbb{P}$. We must find a $G_2$ as in the statement of the lemma. Working in $M[G_1]$, set

$$\mathbb{Q} = \{\, p \in \mathbb{Q}'_2 \mid p \text{ is compatible with every element of } G_1 \,\}.$$

**Claim 2.** *In $M[G_1]$, r.o.$(\mathbb{P})$ can be completely embedded into r.o.$(\mathbb{Q})$.*

*Proof of Claim 2.* Work in $\mathcal{M}[G_1]$. Let $\gamma = \mathrm{card}\big(\mathcal{P}(\mathrm{r.o.}(\mathbb{P}))\big)$. As $\delta_2$ is inaccessible in $M[G_1]$, $\gamma < \delta_2$. We claim that $\mathbb{Q}$ collapses $\gamma$ to be countable. For let $G$ be $\mathcal{M}[G_1]$-generic over $\mathbb{Q}$. By exercise (D4) on page 244 of [Kunen], $G$ is $M$-generic over $\mathbb{Q}'_2$, and $M[G_1][G] = M[G]$. It then follows that $G$ is $M$-generic over $\mathbb{Q}_2$. By chapter 9 of [Ma1], $\delta_2 = (\omega_1)^{\mathcal{M}[G]}$. So $\gamma$ is countable in $M[G] = M[G_1][G]$. In other words, $\mathbb{Q}$ collapses $\gamma$ to be countable. Our conclusion easily follows from this. For let $\dot{H}$ be a $\mathbb{Q}$-name for a generic for r.o.$(\mathbb{P})$. If $b \in \mathrm{r.o.}(\mathbb{P}))$, let $\pi(b)$ be the boolean value in r.o.$(\mathbb{Q})$ that $\check{b} \in \dot{H}$. As in the proof of Lemma 25.5 from [Jech], it is easy to see that $\pi$ is a complete embedding. $\diamondsuit$

It follows that there is a $G_2$ such that $G_2$ is $M[G_1]$-generic over $\mathbb{Q}$, with $H \in M[G_1][G_2]$. Fix such a $G_2$. By exercise (D4) on page 244 of [Kunen], $G_2$ is $M$-generic over $\mathbb{Q}'_2$ and $M[G_2] = M[G_1][G_2]$. It follows that $G_2$ is $M$-generic over $\mathbb{Q}_2$ and $\hat{S} \in G_2$. Thus $G_2 \cap \mathbb{Q}_1 = G_1$ and our lemma is proved. $\square$

We shall use heavily a result due to H. Woodin concerning genericity over $L[\vec{E}]$ models. This result is often referred to as "iterating to make a real generic." We state the result here in the form in which we shall need it. The following is also Theorem 4.3 from [St2].

**Lemma 2.4** (Woodin). *Let $\mathcal{M}$ be a countable, realizable, premouse. Suppose $\gamma < \delta < \mathrm{OR}^{\mathcal{M}}$. Suppose $\mathcal{M} \models$ "$\delta$ is a Woodin cardinal." Let $\mathbb{P} \subseteq J_\gamma^{\mathcal{M}}$ be a partial order in $\mathcal{M}$ and let $G$ be $\mathcal{M}$-generic over $\mathbb{P}$. Then there is a partial order $\mathbb{Q} \subseteq J_\delta^{\mathcal{M}}$, with $\mathbb{Q} \in \mathcal{M}$, such that for any real $w$, there is an iteration tree $\mathcal{T}$ on $\mathcal{M}$ of countable length $\theta + 1$ such that*

(a) $\mathcal{M}_\theta^{\mathcal{T}}$ *is realizable, and*

(b) $D^{\mathcal{T}} = \emptyset$ *so that $i_{0,\theta}^{\mathcal{T}}$ is defined, and*

(c) $\mathrm{crit}(E_\xi^{\mathcal{T}}) > \gamma$ *for all $\xi < \theta$ (so $G$ is $\mathcal{M}_\theta^{\mathcal{T}}$-generic over $\mathbb{P}$), and*

(d) $w$ is $\mathcal{M}_\theta^{\mathcal{T}}[G]$-generic over $i_{0,\theta}^{\mathcal{T}}(\mathbb{Q})$.

For an idea of how to prove this lemma, see the exercises at the end of chapters 6 and 7 in [Ma1].

Our next two results will say that if a premouse $\mathcal{M}$ has $\omega$ Woodin cardinals cofinal in its ordinals, then $\mathcal{M}$ is, in a certain sense, $\Sigma_1(J_\kappa(\mathbb{R}))$-correct. The following is the main technical lemma in this section. Our proof of this lemma is very similar to, and owes its main idea to Steel's proof of Corollary 4.7 in [St2]. In that proof, steel iterates a premouse $\mathcal{M}$, yielding a premouse $\mathcal{M}'$, with the property that every real in $V$ is generic over $\mathcal{M}'$. Steel then takes a generic ultrapower of $\mathcal{M}'$ yielding an embedding $j : \mathcal{M}' \to \text{Ult}$, with the property that $\mathbb{R}^{\text{Ult}} = \mathbb{R}^V$. Now Ult will in general not be wellfounded, but because of the large cardinals in $\mathcal{M}'$ Steel can arrange that the ordinal height of the wellfounded part of Ult is as large as he wants. (Steel assumes that $\mathcal{M}$ has $\omega$ Woodin cardinals, *plus* another extender above the $\omega$ Woodin cardinals.) This allows Steel to reflect any $\Sigma_1$ truth in $L(\mathbb{R})$ down to $\mathcal{M}'$, and hence down to $\mathcal{M}$. In our setting, we are only assuming that $\mathcal{M}$ has $\omega$ Woodin cardinals cofinal in its ordinals. We will therefore not be able to arrange that the ordinal height of the wellfounded part of Ult is as large as we want. Instead, we will only be able to quote abstract admissibility theory to conclude that the ordinal height of the wellfounded part of Ult is at least $\kappa^{\mathbb{R}}$. This will allow us to reflect $\Sigma_1$ truth from $J_\kappa(\mathbb{R})$ down to $\mathcal{M}'$, and hence down to $\mathcal{M}$.

**Lemma 2.5.** *Let $\mathcal{M}$ be a countable, realizable, meek premouse. Suppose that $\delta_1 < \delta_2 < \delta_3 < \cdots$ is an increasing $\omega$-sequence of ordinals such that each $\delta_i$ is a Woodin cardinal of $\mathcal{M}$, and the $\delta_i$ are cofinal in the ordinals of $\mathcal{M}$. Let $\mathbb{P}$ be any partial order in $J_{\delta_1}^{\mathcal{M}}$ and suppose that $H$ is $\mathcal{M}$-generic over $\mathbb{P}$, with $H \in V$. Let $\mathbb{R}^* = \mathbb{R} \cap \mathcal{M}[H]$. Then there is an ordinal $\alpha < \delta_1$ such that whenever $x \in \mathbb{R}^*$ and $\varphi$ is a $\Sigma_1$ formula and $J_\kappa(\mathbb{R}) \models \varphi[x, \mathbb{R}]$, then $J_\alpha(\mathbb{R}^*) \models \varphi[x, \mathbb{R}^*]$.*

*Proof.* Fix an ordinal $\mu < \delta_1$ such that $\mathbb{P} \in J_\mu^{\mathcal{M}}$. Let $\text{Coll}(\omega, 2^{\aleph_0})$ be the partial order which collapses $2^{\aleph_0}$ to be countable. For the rest of the proof we work in a generic extension of $V$ via this partial order. Let $x_1, x_2, x_3, \ldots$ be an $\omega$-enumeration of the reals of $V$.

We are going to iterate $\mathcal{M}$ to yield a new premouse $\mathcal{M}'$ with the property that every real of $V$ is generic over $\mathcal{M}'[H]$. More specifically, we are going to define a sequence of mice: $\mathcal{M}_0, \mathcal{M}_1, \mathcal{M}_2, \ldots$ with $\mathcal{M}_0 = \mathcal{M}$, and a commutative system of embeddings $j_{n,m} : \mathcal{M}_n \to \mathcal{M}_m$, with the property that $\text{crit}(j_{0,n}) > \mu$ so that $H$ is $\mathcal{M}_n$-generic over $\mathbb{P}$, and a sequence of partial orders $\mathbb{Q}_1, \mathbb{Q}_2, \mathbb{Q}_3, \ldots$ such that for each $n \geq 1$, $\mathbb{Q}_n \in \mathcal{M}_n[H]$, and a sequence of generic objects $G_1, G_2, G_3, \ldots$ such that for each $n \geq 1$, $G_n$ is $\mathcal{M}_n[H]$-generic over $\mathbb{Q}_n$. We will arrange that for each $n \geq 1$, $x_n \in \mathcal{M}_n[H][G_n]$. We will maintain inductively that every finite initial segment of our construction lives in $V$. The whole construction will not live in $V$. Then we will define $\mathcal{M}'$ to be the direct limit of the $\mathcal{M}_n$.

Given such a construction, let $\bar{\delta}_{2n} = j_{0,n}(\delta_{2n})$. We will also arrange that:

- $\mathbb{Q}_n \subseteq V_{\bar{\delta}_{2n}}^{\mathcal{M}_n[H]}$, and
- $\text{crit}(j_{n,m}) > \bar{\delta}_{2n}$.

Thus if $m > n$ then we will have that $G_n$ is $\mathcal{M}_m[H]$-generic over $\mathbb{Q}_n$. Finally, we will arrange that if $m > n$ then $G_m \cap \mathbb{Q}_n = G_n$.

We begin by setting $\mathcal{M}_0 = \mathcal{M}$. Now, working in $V$ and using Lemma 2.4 above, we can find a partial order $\hat{\mathbb{P}} \subseteq J_{\delta_1}^{\mathcal{M}_0}$, with $\hat{\mathbb{P}} \in \mathcal{M}_0$, and an iteration tree $\mathcal{T}$ on $\mathcal{M}_0$ of countable length $\theta + 1$ such that

(a) $\mathcal{M}_\theta^\mathcal{T}$ is realizable, and
(b) $D^\mathcal{T} = \emptyset$ so that $i_{0,\theta}^\mathcal{T}$ is defined, and
(c) $\text{crit}(E_\xi^\mathcal{T}) > \mu$ for all $\xi < \theta$ (so $H$ is $\mathcal{M}_\theta^\mathcal{T}$-generic over $\mathbb{P}$), and
(d) $x_1$ is $\mathcal{M}_\theta^\mathcal{T}[H]$-generic over $i_{0,\theta}^\mathcal{T}(\hat{\mathbb{P}})$.

Set $\mathcal{M}_1 = \mathcal{M}_\theta^\mathcal{T}$ and set $j_{0,1} = i_{0,\theta}^\mathcal{T}$. Now let $\bar{\delta}_2 = j_{0,1}(\delta_2)$, and let $\mathbb{Q}_1 = \mathbb{Q}_{\bar{\delta}_2}^{\mathcal{M}_1[H]}$. Finally, again working in $V$ and using (part of the statement of) Lemma 2.3, let $G_1$ be an $\mathcal{M}_1[H]$-generic filter over $\mathbb{Q}_1$ such that $x_1 \in \mathcal{M}_1[H][G_1]$.

Now we describe the second step of our construction. Working in $V$ again, and using Lemma 2.4, we can find a partial order $\hat{\mathbb{P}} \subseteq J_{j_{0,1}(\delta_3)}^{\mathcal{M}_1}$, with $\hat{\mathbb{P}} \in \mathcal{M}_1$, and an iteration tree $\mathcal{T}$ on $\mathcal{M}_1$ of countable length $\theta + 1$ such that

(a) $\mathcal{M}_\theta^{\mathcal{J}}$ is realizable, and

(b) $D^{\mathcal{J}} = \emptyset$ so that $i_{0,\theta}^{\mathcal{J}}$ is defined, and

(c) $\text{crit}(E_\xi^{\mathcal{J}}) > \bar{\delta}_2$ for all $\xi < \theta$ (so $H$ is $\mathcal{M}_\theta^{\mathcal{J}}$-generic over $\mathbb{P}$ and $G_1$ is $\mathcal{M}_\theta^{\mathcal{J}}[H]$-generic over $\mathbb{Q}_1$), and

(d) $x_2$ is $\mathcal{M}_\theta^{\mathcal{J}}[H][G_1]$-generic over $i_{0,\theta}^{\mathcal{J}}(\hat{\mathbb{P}})$.

Set $\mathcal{M}_2 = \mathcal{M}_\theta^{\mathcal{J}}$ and set $j_{1,2} = i_{0,\theta}^{\mathcal{J}}$ and set $j_{0,2} = j_{0,1} \circ j_{1,2}$. Now let $\bar{\delta}_4 = i_{0,2}(\delta_4)$, and let $\mathbb{Q}_2 = \mathbb{Q}_{\bar{\delta}_4}^{\mathcal{M}_2[H]}$. Finally, again working in $V$ and using Lemma 2.3, let $G_2$ be an $\mathcal{M}_2[H]$-generic filter over $\mathbb{Q}_2$ such that $x_2 \in \mathcal{M}_2[H][G_2]$, and $G_2 \cap \mathbb{Q}_1 = G_1$.

The $n$th step of the construction for $n > 2$ is similar to the second step of the construction. It is easy to see how to continue the construction so as to obtain a sequence of mice $\mathcal{M}_0, \mathcal{M}_1, \mathcal{M}_2, \ldots$, a commutative system of embeddings $j_{n,m} : \mathcal{M}_n \to \mathcal{M}_m$, a sequence of partial orders $\mathbb{Q}_1, \mathbb{Q}_2, \mathbb{Q}_3, \ldots$, and a sequence of generic filters $G_1, G_2, G_3, \ldots$ as described above. Letting $\bar{\delta}_{2n} = j_{0,n}(\delta_{2n})$ we have, in summary, the following:

(1) $\text{crit}(j_{0,n}) > \mu$ so $H$ is $\mathcal{M}_n$-generic over $\mathbb{P}$.

(2) $\mathbb{Q}_n = \mathbb{Q}_{\bar{\delta}_{2n}}^{\mathcal{M}_n[H]}$ for $n \geq 1$.

(3) $G_n$ is $\mathcal{M}_n[H]$-generic over $\mathbb{Q}_n$, for $n \geq 1$.

(4) $x_n \in \mathcal{M}_n[H][G_n]$, for $n \geq 1$.

(5) For $1 \leq n < m$, $\text{crit}(j_{n,m}) > \bar{\delta}_{2n}$ so $G_n$ is $\mathcal{M}_m[H]$-generic over $\mathbb{Q}_n$.

(6) For $1 \leq n < m$, $G_m \cap \mathbb{Q}_n = G_n$.

(7) For $n \geq 0$, $\mathcal{M}_n[H][G_n] \in V$.

Let $\mathcal{M}'$ be the direct limit of the $\mathcal{M}_n$ under the $j_{n,m}$. It is obvious that $\mathcal{M}'$ is wellfounded, because every thread in the direct limit system is eventually constant. Let $j : \mathcal{M} \to \mathcal{M}'$ be the direct limit map. The following facts follow easily from (1) through (7) above:

(a) $\bar{\delta}_{2n} = j(\delta_{2n})$.

(b) $H$ is $\mathcal{M}'$-generic over $\mathbb{P}$.

(c) $\mathbb{Q}_n = \mathbb{Q}_{\bar{\delta}_{2n}}^{\mathcal{M}'[H]}$.

(d) $G_n$ is $\mathcal{M}'[H]$-generic over $\mathbb{Q}_n$.

(e) $x_n \in \mathcal{M}'[H][G_n]$.

(f) $\mathbb{R} \cap \mathcal{M}'[H][G_n] = \mathbb{R} \cap \mathcal{M}_n[H][G_n] \in V$, for $n \geq 1$.

It follows from (e) and (f) above that $\bigcup_n \mathbb{R} \cap \mathcal{M}'[H][G_n] = \mathbb{R}^V$.

Our next step is to form the generic ultrapower of $\mathcal{M}'[H]$ via the $G_n$s. Let $\mathcal{P}_n = \text{Ult}(\mathcal{M}'[H], G_n)$. Since $\bar{\delta}_{2n}$ is a Woodin cardinal of $\mathcal{M}'[H]$, by Chapter 9 of [Ma1], $\mathcal{P}_n$ is wellfounded. Also, by by Chapter 9 from [Ma1], $\mathbb{R} \cap \mathcal{P}_n = \mathbb{R} \cap \mathcal{M}'[H][G_n]$. Let $i_n : \mathcal{M}'[H] \to \mathcal{P}_n$ be the ultrapower embedding. Then $i_n$ is a cofinal, $\Sigma_0$-embedding. For $1 \leq n < m$, let $k_{n,m} : \mathcal{P}_n \to \mathcal{P}_m$ be the canonical embedding. That is:

$$k_{n,m}([f]_{G_n}) = [f]_{G_m}.$$

Then the $k_{n,m}$ form a commutative system of embeddings. Also, for each $n < m$, $i_m = k_{n,m} \circ i_n$. Let $\mathcal{P}$ be the direct limit of the $\mathcal{P}_n$ under the $k_{n,m}$ embeddings. Let $k_n : \mathcal{P}_n \to \mathcal{P}$ be the direct limit map. Also let $i : \mathcal{M}'[H] \to \mathcal{P}$ be the induced map. We may think of $\mathcal{P}$ as being the ultrapower of $\mathcal{M}'[H]$ via $\bigcup_n G_n$, and we may think of $i$ as the ultrapower embedding.

Now $\mathcal{P}$ will, in general, not be wellfounded. Let us identify the wellfounded part of $\mathcal{P}$ with the transitive set to which it is isomorphic. It is easy to see that $\mathbb{R}^{\mathcal{P}} = \bigcup_n \mathbb{R} \cap \mathcal{P}_n = \mathbb{R}^V$. It follows from general admissibility theory that the rank of the wellfounded part of $\mathcal{P}$ is at least $\kappa^{\mathbb{R}}$.

Now fix a real $x \in \mathbb{R}^*$, and fix a $\Sigma_1$ formula $\varphi$ such that $J_\kappa(\mathbb{R}) \models \varphi[x, \mathbb{R}]$. Then there is a $\gamma < \kappa^{\mathbb{R}}$ such that $J_\gamma(\mathbb{R}) \models \varphi[x, \mathbb{R}]$. Since $\gamma$ is in the wellfounded part of $\mathcal{P}$, $\mathcal{P}$ satisfies the following $\Sigma_1$ sentence:

> There is an ordinal $\gamma$ such that $J_\gamma(\mathbb{R}) \models \varphi[x, \mathbb{R}]$ and for all $\eta \leq \gamma$, $J_\eta(\mathbb{R})$ is not admissible.

Since $i$ is a cofinal $\Sigma_0$ embedding, $\mathcal{M}'[H]$ also satisfies this sentence. Fix a name $\dot{x}$ in $\mathcal{M}'$ such that $x = \dot{x}_H$. Since $\text{crit}(j) > \mu$ we may pick $\dot{x}$ so that it is fixed by $j$. Let $p \in H$ be such that in $\mathcal{M}'$, $p$ forces the statement:

> There is an ordinal $\gamma$ such that $J_\gamma(\mathbb{R}) \models \varphi[\dot{x}, \mathbb{R}]$ and for all $\eta \leq \gamma$, $J_\eta(\mathbb{R})$ is not admissible.

Again, since $\text{crit}(j) > \mu$, we have that $j(p) = p$. Thus $p$ forces the same statement over $\mathcal{M}$. So in $\mathcal{M}[H]$ we have that There is an ordinal $\gamma$ such that $J_\gamma(\mathbb{R}) \models \varphi[x, \mathbb{R}]$ and for all $\eta \leq \gamma$, $J_\eta(\mathbb{R})$ is not admissible.

Fix such an ordinal $\gamma$ in $\mathcal{M}[H]$. Since $\mathbb{R} \cap \mathcal{M}[H] = \mathbb{R}^*$, we have that $J_\gamma(\mathbb{R}^*) \models \varphi[x, \mathbb{R}^*]$. Since $J_\eta(\mathbb{R}^*)$ is not admissible for all $\eta \leq \gamma$, we have that $\gamma$ is less than, for instance, $\left((2^{\aleph_0})^+\right)^{\mathcal{M}[H]}$. This completes the proof of the lemma. For the ordinal $\alpha$ mentioned in the statement of the lemma we can take, say, $\alpha = (\mu^{++})^{\mathcal{M}}$. □

The previous Lemma gave us an $\alpha$ such that $\Sigma_1(J_\kappa(\mathbb{R}))$ goes down to $J_\alpha(\mathbb{R}^*)$. The next corollary says that $\alpha$ can be chosen so that $\Sigma_1(J_\kappa(\mathbb{R}))$ also goes up.

**Corollary 2.6.** *Under the hypotheses of the previous Lemma, there is an ordinal $\alpha < \delta_1$ such that whenever $x \in \mathbb{R}^*$ and $\varphi$ is a $\Sigma_1$ formula then $J_\kappa(\mathbb{R}) \models \varphi[x, \mathbb{R}]$ iff $J_\alpha(\mathbb{R}^*) \models \varphi[x, \mathbb{R}^*]$.*

*Proof.* Let $\alpha_0$ be given by the previous lemma. If $\alpha_0$ does not satisfy our corollary, then there is some $\Sigma_1$ formula $\varphi$, and some $x \in \mathbb{R}^*$ such that $J_{\alpha_0}(\mathbb{R}^*) \models \varphi[x, \mathbb{R}^*]$, and $J_\kappa(\mathbb{R}) \not\models \varphi[x, \mathbb{R}]$. Let $\gamma \leq \alpha_0$ be least such that there is some $\Sigma_1$ formula $\varphi$, and some $x \in \mathbb{R}^*$ such that $J_\gamma(\mathbb{R}^*) \models \varphi[x, \mathbb{R}^*]$, and $J_\kappa(\mathbb{R}) \not\models \varphi[x, \mathbb{R}]$. Then $\gamma$ is a successor ordinal. Let $\alpha_1$ be such that $\gamma = \alpha_1 + 1$. We claim that $\alpha_1$ witnesses that our corollary is true. So let $\varphi$ by any $\Sigma_1$ formula and let $x \in \mathbb{R}^*$. If $J_{\alpha_1}(\mathbb{R}^*) \models \varphi[x, \mathbb{R}^*]$, then by definition of $\gamma$, $J_\kappa(\mathbb{R}) \models \varphi[x, \mathbb{R}]$. Conversely, suppose that $J_\kappa(\mathbb{R}) \models \varphi[x, \mathbb{R}]$. Fix $y \in \mathbb{R}^*$, and fix a $\Sigma_1$ formula $\psi$ so that $J_\gamma(\mathbb{R}^*) \models \psi[y, \mathbb{R}^*]$ but $J_\kappa(\mathbb{R}) \not\models \psi[y, \mathbb{R}]$. Then $J_\kappa(\mathbb{R})$ satisfies the following $\Sigma_1$ statement about $\langle x, y \rangle$:

There is an ordinal $\beta$ such that $J_\beta(\mathbb{R}) \models \varphi[x, \mathbb{R}]$ and $J_\beta(\mathbb{R}) \not\models \psi[y, \mathbb{R}]$.

Thus $J_{\alpha_0}(\mathbb{R}^*)$ satisfies the same statement about $\langle x, y \rangle$. So there is an ordinal $\beta < \alpha_0$ such that $J_\beta(\mathbb{R}^*) \models \varphi[x, \mathbb{R}^*]$, and $J_\beta(\mathbb{R}^*) \not\models \psi[y, \mathbb{R}^*]$. Since $J_\gamma(\mathbb{R}^*) \models \psi[y, \mathbb{R}^*]$, we have that $\beta \leq \alpha_1$. Thus $J_{\alpha_1}(\mathbb{R}^*) \models \varphi[x, \mathbb{R}^*]$. □

**Corollary 2.7.** *Let $\mathcal{M}$ be a countable, realizable, meek premouse. Suppose that $\delta_1 < \delta_2 < \delta_3 < \cdots$ is an increasing $\omega$-sequence of ordinals such that each $\delta_i$ is a Woodin cardinal of $\mathcal{M}$, and the $\delta_i$ are cofinal in the ordinals of $\mathcal{M}$. Let $x$ be a real and suppose that for some $\beta < \kappa^{\mathbb{R}}$, $x$ is ordinal definable over $J_\beta(\mathbb{R})$. Then $x \in \mathcal{M}$.*

*Proof.* By Lemma 2.4 there is a partial order $\mathbb{Q} \subseteq J_{\delta_1}^{\mathcal{M}}$, with $\mathbb{Q} \in \mathcal{M}$ and there is an iteration tree $\mathcal{T}$ on $\mathcal{M}$ of countable length $\theta + 1$ such that

(a) $\mathcal{M}_\theta^\mathcal{T}$ is realizable, and

(b) $D^\mathcal{T} = \emptyset$ so that $i_{0,\theta}^\mathcal{T}$ is defined, and

(c) $x$ is $\mathcal{M}_\theta^\mathcal{T}$-generic over $i_{0,\theta}^\mathcal{T}(\mathbb{Q})$.

Let $\mathcal{M}' = \mathcal{M}_\theta^\mathcal{T}$ and let $\delta_1' = i_{0,\theta}^\mathcal{T}(\delta_1)$. It suffices to show that $x \in \mathcal{M}'$. Let $\mathbb{P}$ be be the partial order in $\mathcal{M}'$ for collapsing $\delta_1'$ to be countable. Then there is a filter $H$ which is $\mathcal{M}'$-generic over $\mathbb{P}$ with $x \in \mathcal{M}'[H]$. Let $\mathbb{R}^* = \mathbb{R} \cap \mathcal{M}'[H]$. Let us apply the previous corollary to $\mathcal{M}'$ and $H$. This gives us an ordinal $\alpha \in \mathcal{M}'$ such that $J_\alpha(\mathbb{R}^*)$ satisfies the following $\Sigma_1$ statement about $x$:

There is a $\beta$ such that $x$ is ordinal definable over $J_\beta(\mathbb{R})$.

Thus there is a $\beta$ in $\mathcal{M}'$ so that $x$ is ordinal definable over $J_\beta(\mathbb{R}^*)$. This means that $x$ is ordinal definable in $\mathcal{M}'[H]$. As $\mathbb{P}$ is a homogeneous forcing, $x \in \mathcal{M}'$. (Actually, since $\mathcal{M}'$ is not a model of ZFC, we are really using that $x$ is ordinal definable in $J_\gamma^{\mathcal{M}'}[H]$, for some $\gamma$ such that $J_\gamma^{\mathcal{M}'} \models \text{ZFC}$.) □

**Theorem 2.8.** *Assume that there are $\omega$ Woodin cardinals. Let $\mathcal{M}$ be a countable, realizable premouse, such that $\mathcal{M}$ is not fairly small, but every proper initial segment of $\mathcal{M}$ is fairly small. Then the reals of $\mathcal{M}$ are exactly equal to the set of reals $x$ such that $x$ is ordinal definable over $J_\beta(\mathbb{R})$, for some $\beta < \kappa^\mathbb{R}$.*

*Proof.* This follows immediately from Theorem 1.6 on page 10, and Corollary 2.7 above. □

**Remark 2.9.** Assume that there are $\omega$ Woodin cardinals. It is shown in [St2] that this implies that there is a proper class premouse $L[\vec{E}]$ which satisfies ZFC +"There are $\omega$ Woodin cardinals." Let $\mathcal{M}$ be the $\trianglelefteq$-least initial segment of $L[\vec{E}]$ which is fairly big. It is easy to see that $\mathcal{M}$ projects to $\omega$. $\mathcal{M}$ is thus the unique fully iterable premouse which is sound, projects to $\omega$, has $\omega$ Woodin cardinals cofinal in its ordinals, and has the property that every proper initial segment is fairly small. We may think of $\mathcal{M}$ as the canonical, minimal inner model for the theory: ZFC − Replacement + "There exists $\omega$ Woodin cardinals cofinal in the ordinals." The previous theorem says that $A_{\text{HYP}} = \mathbb{R} \cap \mathcal{M}$.

## 3. $\Sigma_n^*$ Correctness: A Proof of Martin's Theorem

In the previous section we showed that $A_{\text{HYP}} = \mathbb{R} \cap \mathcal{M}$, where $\mathcal{M}$ is the canonical, minimal inner model for the theory ZFC − Replacement + "There are $\omega$ Woodin cardinals cofinal in the ordinals." In this section we give one application of this characterization of $A_{\text{HYP}}$. We show that every real $x$ which is definable over $J_\kappa(\mathbb{R})$ is in $A_{\text{HYP}}$. We do this by showing that $x \in \mathcal{M}$. Thus we are giving a purely inner-model-theoretic proof of Martin's theorem that every $\Sigma_n^*$ real is in the largest countable inductive set. The reason that the real $x$ is in $\mathcal{M}$ is that for each $n \in \omega$, $\mathcal{M}$ can compute $\Sigma_n(J_\kappa(\mathbb{R}))$ truth. This correctness result for $\mathcal{M}$ is the main content of this section.

Let $\mathcal{M}$ be a countable, realizable premouse such that $\mathcal{M}$ is fairly big, but every proper initial segment of $\mathcal{M}$ is fairly small. Corollary 2.6 tells us that there is an ordinal $\alpha$ in $\mathcal{M}$ so that $(J_\alpha(\mathbb{R}))^{\mathcal{M}}$ agrees with $J_\kappa(\mathbb{R})$ on $\Sigma_1$ facts about reals in $\mathcal{M}$. Consider the following question. Is it possible to identify $\alpha$ while working in $\mathcal{M}$? The next lemma tells us that the answer is yes.

**Lemma 3.1.** *Assume that there exists $\omega$ Woodin cardinals in $V$. Let $\mathcal{M}$ be a countable, realizable premouse such that $\mathcal{M}$ is fairly big, but every proper initial segment of $\mathcal{M}$ is fairly small. Let $\gamma$ be any cardinal of $\mathcal{M}$ and suppose that $G$ is $\mathcal{M}$-generic over $\mathrm{Coll}(\omega, \gamma)$, with $G \in V$. Let $\mathbb{R}^* = \mathbb{R} \cap \mathcal{M}[G]$. Suppose that $\alpha_0 \in \mathcal{M}$ is an ordinal such that whenever $x \in \mathbb{R}^*$ and $\varphi$ is a $\Sigma_1$ formula then $J_\kappa(\mathbb{R}) \models \varphi[x, \mathbb{R}]$ iff $J_{\alpha_0}(\mathbb{R}^*) \models \varphi[x, \mathbb{R}^*]$. Then $\alpha_0$ is the least $\alpha \in \mathcal{M}$ such that there is a wellorder of $\mathbb{R}^*$ which is definable over $J_\alpha(\mathbb{R}^*)$.*

*Proof.* First let us see that there is a wellorder of $\mathbb{R}^*$ which is definable over $J_{\alpha_0}(\mathbb{R}^*)$. Every element of $\mathbb{R}^*$ has a $G$-name in $J_\gamma^{\mathcal{M}}$. Furthermore $J_\gamma^{\mathcal{M}}$ and $G$ are in $J_{\alpha_0}(\mathbb{R}^*)$ as they are coded by elements of $\mathbb{R}^*$. So it suffices to see that there is a wellorder of $J_\gamma^{\mathcal{M}}$ which is definable over $J_{\alpha_0}(\mathbb{R}^*)$. Let $<$ be the order of construction in $\mathcal{M}$. We claim that $< \restriction J_\gamma^{\mathcal{M}}$ is definable over $J_{\alpha_0}(\mathbb{R}^*)$. Let $\psi$ be a $\Pi_1$ formula such that for all $w \in \mathbb{R}$, $J_\kappa(\mathbb{R}) \models \psi[w, \mathbb{R}]$ iff $w$ codes a countable premouse $\mathcal{N}_w$, and $\mathcal{N}_w$ is $\partial^{\mathbb{R}}$-closed iterable. It follows from our comparison theorem, Lemma 1.5 on page 9, that for $x, y \in J_\gamma^{\mathcal{M}}$, $x < y$ iff $J_{\alpha_0}(\mathbb{R}^*)$ satisfies the following statement:

> There is a $w \in \mathbb{R}^*$ such that $\psi[w, \mathbb{R}^*]$, and $x, y \in \mathcal{N}_w$ and $x$ is less than $y$ in the order of construction of $\mathcal{N}_w$.

So we have shown that there is a wellorder of $\mathbb{R}^*$ which is definable over $J_{\alpha_0}(\mathbb{R}^*)$.

To see that $\alpha_0$ is least with this property, we must show that there is no wellorder of $\mathbb{R}^*$ which is *in* $J_{\alpha_0}(\mathbb{R}^*)$. By $\Sigma_1$-correctness, it suffices to see that there is no wellorder of $\mathbb{R}$ in $J_\kappa(\mathbb{R})$. But this follows from our assumption that there are $\omega$ Woodin cardinals in $V$. In fact our large cardinal hypothesis implies that every game in $J_\kappa(\mathbb{R})$ (and more) is determined. [proof: Using Woodin's stationary tower forcing as in the end of the proof of our comparison lemma, Lemma 1.5, we get a fully elementary generic embedding $j : J_\kappa(\mathbb{R}) \to J_\lambda(\mathbb{R}')$, where $\lambda$ is some ordinal, and where $\mathbb{R}'$ is the set of reals in a symmetric collapse of $V$ up to the sup of the $\omega$ Woodin cardinals. Also using stationary tower forcing, Woodin has shown that $L(\mathbb{R}') \models$ AD. (See chapter 9 of [Ma1].) In particular $J_\lambda(\mathbb{R}') \models$ AD, and so $J_\kappa(\mathbb{R}) \models$ AD.] $\square$

The following is the main technical result of this section. The lemma says that if $\mathcal{M}$ has $\omega$ Woodin cardinals cofinal in its ordinals, then $\mathcal{M}$ can compute $\Sigma_n(J_\kappa(\mathbb{R}))$ truth, for every $n$. Our proof is by induction on $n$. In order to carry out this inductive proof we need to make the inductive hypothesis *uniform* in $\mathcal{M}$. Unfortunately, this slightly complicates the statement of the lemma.

**Lemma 3.2.** *Assume that there exists $\omega$ Woodin cardinals in $V$. Let $n \geq 1$, and let $\varphi$ be a $\Sigma_n$ formula. Then there is another formula $\psi = \psi_\varphi$ such that whenever $\mathcal{M}$ is a countable, realizable premouse, and $\mathcal{M}$ is fairly big but every proper initial segment of $\mathcal{M}$ is fairly small, and $\gamma$ is a cardinal of $\mathcal{M}$, and $G$ is $\mathcal{M}$-generic over $\mathrm{Coll}(\omega, \gamma)$, with $G \in V$, and $\delta_1, \delta_2, \ldots \delta_n$ are Woodin cardinals of $\mathcal{M}$ with $\gamma < \delta_1 < \delta_2 < \cdots < \delta_n$, and $\mathbb{R}^* = \mathbb{R} \cap \mathcal{M}[G]$, and $x \in \mathbb{R}^*$, then $J_\kappa(\mathbb{R}) \models \varphi[x, \mathbb{R}]$ iff $J_{\delta_n}^{\mathcal{M}}[G] \models \psi[x, \delta_1, \ldots, \delta_{n-1}]$.*

*Proof.* By induction on $n$. First let $n = 1$ and let $\varphi$ be a $\Sigma_1$ formula. Then we can take $\psi = \psi_\varphi$ to be a formula such that, given any $\mathcal{M}$, $\gamma$, $G$ and $\delta_1$ as above, and letting $\mathbb{R}^* = \mathbb{R} \cap \mathcal{M}[G]$, and letting $x \in \mathbb{R}^*$, we have that $J_{\delta_1}^{\mathcal{M}}[G] \models \psi[x]$ iff $(\exists \alpha < \delta_1)$ s.t. there is a wellorder of $\mathbb{R}^*$ definable over $J_\alpha(\mathbb{R}^*)$, and letting $\alpha_0$ be the least such $\alpha$ we have that $J_\alpha(\mathbb{R}^*) \models \varphi[x]$. Clearly there is such a formula $\psi$. Furthermore such a $\psi$ works by Corollary 2.6 and Lemma 3.1.

Next let $n > 1$ and let $\varphi$ be a $\Sigma_n$ formula. Let $\theta$ be a $\Sigma_{n-1}$ formula so that $\varphi(v_1, v_2) \leftrightarrow (\exists v_0) \neg \theta(v_0, v_1, v_2)$. Let $\psi' = \psi_\theta$. ($\psi_\theta$ exists by induction.) Then we can take $\psi = \psi_\varphi$ to be a formula such that, given any $\mathcal{M}$, $\gamma$, $G$ and $\delta_1, \ldots \delta_n$ as above, and letting $\mathbb{R}^* = \mathbb{R} \cap \mathcal{M}[G]$, and letting $x \in \mathbb{R}^*$, we have that $J_{\delta_n}^{\mathcal{M}}[G] \models \psi[x, \delta_1, \ldots \delta_{n-1}]$ iff it is forced over $J_{\delta_n}^{\mathcal{M}}[G]$ that in a generic extension via $\text{Coll}(\omega, \delta_1)$, there is a real $y$ such that $\neg \psi'[y, x, \delta_2, \ldots, \delta_{n-1}]$. Clearly there is such a formula $\psi$. Let us see that such a $\psi$ works.

Fix $\mathcal{M}$, $\gamma$, $G$ and $\delta_1, \ldots, \delta_n$ as above. Let $\mathbb{R}^* = \mathbb{R} \cap \mathcal{M}[G]$, and let $x \in \mathbb{R}^*$. First suppose that $J_{\delta_n}^{\mathcal{M}}[G] \models \psi[x, \delta_1, \ldots \delta_{n-1}]$. Let $H$ be $\mathcal{M}[G]$-generic over $\text{Coll}(\omega, \delta_1)$, and let $\mathbb{R}^{**} = \mathbb{R} \cap \mathcal{M}[G][H]$. By definition of $\psi$, there is a $y \in \mathbb{R}^{**}$ such that $J_{\delta_n}^{\mathcal{M}}[G][H] \models \neg \psi'[x, y, \delta_2, \ldots, \delta_{n-1}]$. Since $\psi' = \psi_\theta$, we have by induction that $J_\kappa(\mathbb{R}) \models \neg \theta[y, x, \mathbb{R}]$. Thus $J_\kappa(\mathbb{R}) \models \varphi[x, \mathbb{R}]$.

Conversely, suppose that $J_\kappa(\mathbb{R}) \models \varphi[x, \mathbb{R}]$. Then there is a real $y$ such that $J_\kappa(\mathbb{R}) \models \neg \theta[y, x, \mathbb{R}]$. Fix such a $y$. By Lemma 2.4, there is is a partial order $\mathbb{Q} \subseteq J_{\delta_1}^{\mathcal{M}}$, with $\mathbb{Q} \in \mathcal{M}$, and there is an iteration tree $\mathcal{T}$ on $\mathcal{M}$ of countable length $\mu + 1$ such that

(a) $\mathcal{M}_\mu^{\mathcal{T}}$ is realizable, and

(b) $D^{\mathcal{T}} = \emptyset$ so that $i_{0,\mu}^{\mathcal{T}}$ is defined, and

(c) $\text{crit}(E_\xi^{\mathcal{T}}) > \gamma$ for all $\xi < \mu$ (so $G$ is $\mathcal{M}_\mu^{\mathcal{T}}$-generic over $\text{Coll}(\omega, \gamma)$), and

(d) $y$ is $\mathcal{M}_\mu^{\mathcal{T}}[G]$-generic over $i_{0,\mu}^{\mathcal{T}}(\mathbb{Q})$.

Let $\mathcal{M}' = \mathcal{M}_\mu^{\mathcal{T}}$, let $i = i_{0,\mu}^{\mathcal{T}}$, and let $\delta_k' = i(\delta_k)$ for $k \geq 1$. Let $H$ be $\mathcal{M}'[G]$-generic over $\text{Coll}(\omega, \delta_1')$ with $y \in \mathcal{M}'[G][H]$. By induction, since $\psi' = \psi_\theta$, $J_{\delta_n'}^{\mathcal{M}'}[G][H] \models \neg \psi'[y, x, \delta_2', \ldots, \delta_{n-1}']$. Since the collapse forcing is homogeneous, it is forced over $J_{\delta_n'}^{\mathcal{M}'}[G]$ that in a generic extension via $\text{Coll}(\omega, \delta_1')$, there is a real $y$ such that $\neg \psi'[y, x, \delta_2', \ldots, \delta_{n-1}']$. Thus $J_{\delta_n'}^{\mathcal{M}'}[G] \models \psi[x, \delta_1', \ldots, \delta_{n-1}']$. Since $\text{crit}(i) > \gamma$, it follows that $J_{\delta_n}^{\mathcal{M}}[G] \models \psi[x, \delta_1, \ldots, \delta_{n-1}]$. □

**Theorem 3.3.** *Assume there exists $\omega$ Woodin cardinals in $V$. Let $\mathcal{M}$ be a countable, realizable premouse such that $\mathcal{M}$ is fairly big, but every proper initial segment of $\mathcal{M}$ is fairly small. Let $x$ be a real which is definable over $J_\kappa(\mathbb{R})$. Then $x \in \mathcal{M}$.*

*Proof.* Let $n \geq 1$ be such that, as a subset of $\omega$, $x$ is $\Sigma_n$ definable over $J_\kappa(\mathbb{R})$. Then by the previous lemma, $x$ is definable over $J_{\delta_n}^{\mathcal{M}}$ from $\delta_1, \ldots, \delta_{n-1}$, where $\delta_1, \cdots, \delta_n$ are the first $n$ Woodin cardinals of $\mathcal{M}$. $\square$

**Corollary 3.4.** *Every $\Sigma_n^*$ real is in the largest countable inductive set.*

Department of Mathematics, Florida International University, Miami, FL 33199

*E-mail address*: rudomine@fiu.edu